\documentclass[a4paper,8pt]{amsart}

\usepackage{amsmath, amsfonts, amsthm, amssymb}
\input xypic

\newtheorem{theorem}{Theorem}
\newtheorem{proposition}[theorem]{Proposition}
\newtheorem{lemma}[theorem]{Lemma}

\theoremstyle{definition}

\newcounter{bean}

\newcommand{\namedright}[3]{\ensuremath{#1\stackrel{#2}
 {\longrightarrow}#3}}
\newcommand{\nameddright}[5]{\ensuremath{#1\stackrel{#2}
 {\longrightarrow}#3\stackrel{#4}{\longrightarrow}#5}}

\newcommand{\qqed}{\hfill\Box}

\begin{document}


\title{The homotopy type of the complement of the codimension-two
       coordinate subspace arrangement}
\author{Jelena Grbi\'{c} and Stephen Theriault}
\address{Department of Mathematical Sciences,
         University of Aberdeen, Aberdeen AB24 3UE, United Kingdom}
\email{jelena@maths.abdn.ac.uk}
\email{s.theriault@maths.abdn.ac.uk}

\renewcommand{\thefootnote}{}


\maketitle


A complex coordinate subspace of $\mathbb{C}^{n}$ is given by
\[L_{\sigma}=\{ (z_{1},\ldots,z_{n})\in \mathbb{C}^{n} |\, z_{i_{1}}=
       \cdots=z_{i_{k}}=0\}\]
where $\sigma =\{i_{1},\ldots,i_{k}\}$ is a subset of $[m]$. For
each simplicial complex $K$ on the set $[m]$ we associate the
complex coordinate subspace arrangement
$\mathcal{C}\mathcal{A}(K)=\{L_{\sigma} |\, \sigma\not\in K\}$ and
its complement $U(K)=\mathbb{C}^{n} \backslash
\bigcup_{\sigma\not\in K} L_{\sigma}$. On the other hand, to $K$
we can associate the Davis -Januszkiewicz space
$DJ(K)=\bigcup_{\sigma\in K}BT_\sigma\subset BT^{n}$, where
$BT^{n}$ is the classifying space of $n$-dimensional torus, that
is, the product of $n$ copies of infinite-dimensional projective
space $\mathbb{C}P^{\infty}$, and ${BT_{\sigma}:=\{ (x_1,\ldots
,x_n)\in BT^n | x_i=\ast}$ where $i \not\in \sigma \}$. Let
$\mathcal{Z}_{K}$ be the fibre of
\(\namedright{DJ(K)}{}{BT^{n}}\). By~\cite[8.9]{BP}, there is an
equivariant deformation retraction
\(\namedright{U(K)}{}{\mathcal{Z}_{K}}\), and the integral
cohomology of $\mathcal{Z}_K$ has been calculated in \cite[7.6 and
7.7]{BP}.

\begin{theorem}
   \label{arrangements}
   The complement of the codimension-two coordinate subspace
   arrangement in $\mathbb{C}^{n}$ has the homotopy type of the wedge
   of spheres
   \[\bigvee_{k=2}^{n} (k-1)\binom{n}{k}S^{k+1}.\]
\end{theorem}

\begin{proof}
Let $K$ be a disjoint union of $n$ vertices. Then $DJ(K)$ is the
wedge of $n$ copies of $\mathbb{C}P^{\infty}$ and $U(K)$ is the
complement of the set of all codimension-two coordinates subspaces
$z_{i}=z_{j}=0$ for $1\leq i<j\leq n$ in $\mathbb{C}^{n}$.
Therefore to prove the theorem we have to determine the homotopy
fibre of the inclusion
\(\namedright{\bigvee_{t=1}^n\mathbb{C}P^{\infty}}{}{\prod^n_{t=1}\mathbb{C}P^\infty}\).
This is done by applying Proposition~\ref{fiberdecomp} to the case
$X_{1}=\cdots=X_{n}=\mathbb{C}P^{\infty}$ and noting that
$\Omega\mathbb{C}P^{\infty}\simeq S^{1}$.
\end{proof}

It should be emphasized that Theorem~\ref{arrangements} holds
without suspending. Previously, decompositions were known only
after some number of suspensions, the best of which was by
Schaper~\cite{S} who required one suspension. To finish the proof
of Theorem~\ref{arrangements} it remains to prove
Proposition~\ref{fiberdecomp}. This was originally proved by
Porter~\cite{P} by examining subspaces of contractible spaces. We
present an accelerated proof based on the Cube Lemma.

We work in the category of based, connected topological spaces and
continuous maps. Let~$\ast$ denote the basepoint. For spaces
$X,Y$, let $X\rtimes Y=(X\times Y)/(\ast\times Y)$, $X\wedge
Y=(X\rtimes Y)/(X\times\ast)$, and $X\ast Y=\Sigma X\wedge Y$.
Denote the identity map on $X$ by $X$. Denote the map which sends
all points to the basepoint by $\ast$.

\begin{lemma}
   \label{podecomp}
   Let $A,B$, and $C$ be spaces. Define $Q$ as the homotopy
   pushout of the map \(\namedright{A\times B}{\ast\times B}{C\times B}\) and
the projection \( \namedright{A\times B}{\pi_{1}}{A}\). Then
$Q\simeq (A\ast B)\vee (C\rtimes B)$.
\end{lemma}

\begin{proof}
Consider the diagram of iterated homotopy pushouts
\[\diagram
     A\times B\rto^-{\pi_{2}}\dto^{\pi_{1}} & B\rto^-{i_{2}}\dto^{\ast}
        & C\times B\dto^{s} \\
     A\rto^-{\ast} & A\ast B\rto^-{t} & \overline{Q}
  \enddiagram\]
where $\pi_{2},i_{2}$ are the projection and inclusion
respectively. Here, it is well known that the left square is a
homotopy pushout, and the right homotopy pushout
defines~$\overline{Q}$. Note that
$i_{2}\circ\pi_{2}\simeq\ast\times B$. The outer rectangle in an
iterated homotopy pushout diagram is itself a homotopy pushout, so
$\overline{Q}\simeq Q$. The right pushout then shows that the
homotopy cofibre of \(\namedright{C\times B}{}{Q}\) is $\Sigma
B\vee (A\ast B)$. Thus~$t$ has a left homotopy inverse. Further,
$s\circ i_{2}\simeq\ast$ so pinching out $B$ in the right pushout
gives a homotopy cofibration \(\nameddright{C\rtimes
B}{}{Q}{r}{A\ast B}\) with $r\circ t$ homotopic to the identity
map.
\end{proof}

\begin{lemma}
   \label{prodsusp}
   Let $Y_{1},\ldots, Y_{n}$ be spaces. Then there is a
   homotopy equivalence
   \[\Sigma (Y_{1}\times\cdots\times Y_{n})\simeq
      \bigvee_{k=1}^{n}\left(\bigvee_{1\leq i_{1}<\cdots <i_{k}\leq n}
      \Sigma Y_{i_{1}}\wedge\cdots\wedge Y_{i_{k}}\right).\]
\end{lemma}

\begin{proof}
Induct on the decomposition $\Sigma (A\times B)=\Sigma A\vee\Sigma
B\vee (\Sigma A\wedge B)$.
\end{proof}

The following was proved by Mather~\cite{M} and is known as the
Cube Lemma.

\begin{lemma}
   \label{cube}
   Suppose there is a diagram of spaces and maps
   \[\diagram
      E\rrto\drto\ddto & & F\dline\drto & \\
      & G\rrto\ddto & \dto & H\ddto \\
      A\rline\drto & \rto & B\drto & \\
      & C\rrto & & D
   \enddiagram\]
   where the bottom face is a homotopy pushout and the four sides are
   obtained by pulling back with
   \(\namedright{H}{}{D}\).
   Then the top face is a homotopy pushout.~$\qqed$
\end{lemma}

\begin{proposition}
   \label{fiberdecomp}
   Let $X_{1},\ldots,X_{n}$ be spaces. Consider the homotopy fibration
   \[\nameddright{F_{n}}{}{X_{1}\vee\cdots\vee X_{n}}{}
        {X_{1}\times\cdots\times X_{n}}\]
   obtained by including the wedge into the product. Then there is a
   homotopy decomposition
   \[F_{n}\simeq\bigvee_{k=2}^{n}\left(
     \bigvee_{1\leq i_{1}<\cdots <i_{k}\leq n} (k-1)
     (\Sigma\Omega X_{i_{1}}\wedge\cdots\wedge\Omega X_{i_{k}})\right).\]
\end{proposition}

\begin{proof}
We induct on $n$. When $n=2$ it is well known that
$F_{2}\simeq\Sigma\Omega X_{1}\wedge\Omega X_{2}$. Let $n\geq 3$
and assume the Proposition holds for $F_{n-1}$. Let
$M_{k}=X_{1}\vee\cdots\vee X_{k}$ and
$N_{k}=X_{1}\times\cdots\times X_{k}$. Observe that $M_{n}$ is the
pushout of $M_{n-1}$ and $X_{n}$ over a point. Composing each
vertex of the pushout into $N_{n}$ we obtain homotopy fibrations
\(\nameddright{\Omega N_{n}}{}{\ast}{}{N_{n}}\),
\(\nameddright{\Omega N_{n-1}}{}{X_{n}}{}{N_{n}}\),
\(\nameddright{F_{n-1}\times\Omega X_{n}}{}{M_{n-1}}{}{N_{n}}\),
and \(\nameddright{F_{n}}{}{M_{n}}{}{N_{n}}\). Write $N_{n}$ as
$N_{n-1}\times X_{n}$. Then Lemma~\ref{cube} implies that there is
a homotopy pushout
\[\diagram
      \Omega N_{n-1}\times\Omega X_{n}\rto^{h}\dto^{g}
           & F_{n-1}\times\Omega X_{n}\dto \\
      \Omega N_{n-1}\rto & F_{n}
  \enddiagram\]
where $g$ is easily identified as the projection and $h$ is the
connecting map for the homotopy fibration
\(\nameddright{F_{n-1}\times\Omega X_{n}}{}{M_{n-1}\times
\ast}{}{N_{n-1}\times X_n}\). So \(
h\simeq\partial_{n-1}\times\Omega X_n\) where \(\partial_{n-1}\)
is the connecting map of the fibration \(
\nameddright{F_{n-1}}{}{M_{n-1}}{}{N_{n-1}}\). But
\(\partial_{n-1}\simeq\ast\) as \(\namedright{\Omega
M_{n-1}}{}{\Omega N_{n-1}}\) has a right homotopy inverse. Thus
$h\simeq\ast\times\Omega X_{n}$. By Lemma~\ref{podecomp},
$F_{n}\simeq(\Omega N_{n-1}\ast\Omega X_{n})\vee
(F_{n-1}\rtimes\Omega X_{n})$. Since $F_{n-1}$ is a suspension,
$F_{n-1}\rtimes\Omega X_{n}\simeq F_{n-1}\vee (F_{n-1}\wedge\Omega
X_{n})$. Combining the decomposition of $\Sigma\Omega
N_{n}\simeq\Sigma (\Omega X_{1}\times\cdots\times\Omega X_{n})$ in
Lemma~\ref{prodsusp} with the inductive decomposition of $F_{n-1}$
and collecting like terms, the asserted wedge decomposition of
$F_{n}$ follows.
\end{proof}

The authors wish to thank Taras Panov and Nigel Ray for
encouraging us to consider this problem.

\bibliographystyle{amsalpha}

\end{document}